\documentclass{article}
\usepackage{amsmath}
\usepackage{amssymb}
\newtheorem{theorem}{Theorem}[section]
\newtheorem{lemma}[theorem]{Lemma}

\newtheorem{remark}{Remark}[section]
\title{Uniqueness of the critical probability for percolation in the two dimensional
 Sierpi\'nski carpet lattice}
\author{Yasunari Higuchi$^{~1,}$\thanks{Research supported in part
by Grant in Aid for Scientific Research No. 15340032} \and
 Xian-Yuan Wu$^{~2,}$\thanks{Research supported in part
by Natural Science Foundation of China in grant No. 10301023} }
\date{}
\begin{document}
\maketitle
\begin{center}
\begin{minipage}{11.5cm}
{\small \noindent\hskip-3mm $^{~1}$Department of Mathematics,
Faculty of Science, Kobe University, Rokko, Kobe 657-8501, Japan.
E-mail: \texttt{higuchi@math.kobe-u.ac.jp}

\vskip 2mm

\noindent\hskip-3mm$^{~2}$Department of Mathematics, Capital Normal
University, 100037, Beijing, China. E-mail:
\texttt{wuxy@mail.cnu.edu.cn}}
\end{minipage}
\end{center}
\noindent {\bf Abstract}\quad We prove that the critical probability
for the Sierpi\'nski carpet lattice in two dimensions is uniquely
determined. The transition is sharp. This extends the Kumagai's
result \cite{Kumagai} to the original Sierpi\'nski carpet
 lattice.

\section{Introduction}
We consider bond percolation problem on the original two dimensional
Sierpi\'nski carpet
 lattice.
Let $T$ be given by $\{ 0,1,2\}^2\setminus \{ (1,1)\} $, and for
each $(i,j)\in T$
 we put
\begin{equation}\label{eq1.1}
\psi_{i,j}(x,y) = 3^{-1}(i,j) + 3^{-1}(x_1,x_2) \qquad \hbox{ for }
(x_1,x_2)\in [0,1]^2
 .
\end{equation}
The Sierpi\'nski carpet is the closed subset $K^T$ of $[0,1]^2$
which satisfies
\begin{equation}\label{eq1.2}
K^T = \bigcup_{(i,j)\in T} \psi_{i,j}(K^T).
\end{equation}
This set is a decreasing limit of
\begin{equation}\label{eq1.3}
K^T_n = \bigcup_{(i_n,j_n)\in T} \cdots \bigcup_{(i_1,j_1)\in T}
\psi_{i_n,j_n}\circ
 \cdots
          \circ \psi_{i_1,j_1}([0,1]^2).
\end{equation}
The pre-Sierpi\'nski carpet lattice $G^T$ is the subgraph of
$\mathbb{Z}^2$, which
 is the
increasing limit of
\begin{equation}\label{eq1.4}
G^T_n = \mathbb{Z}^2\cap 3^n K^T_n.
\end{equation}
To be more precise, (\ref{eq1.4}) determines the vertex set of
$G_n^T$, and the edge set of $G_n^T$ is given in an obvious way;
\begin{equation}\label{eq1.5}
{\mathcal E}(G_n^T)=\left\{ \{ x,y\} \subset G_n^T \  ;  |x-y| =1
\right\} ,
\end{equation}
where $|x-y|$ denotes the Euclidean distance of $x$ and $y$. We will
abuse the notation $G_n^T$ and $G^T$ both for the graphs defined
above and
 their vertex sets.
This will not cause a problem. We define the full Sierpi\'nski
carpet lattice $S^T$, which we sometimes call the original
Sierpi\'nski carpet lattice, by
\begin{equation}\label{eq1.6}
S^T = G^T \cup \Phi_1(G^T) \cup \Phi_2(G^T) \cup
\Phi_1\circ\Phi_2(G^T),
\end{equation}
where $\Phi_j$ is the reflection with respect to the $x_j$-axis for
$j=1,2$, respectively.

Each edge $e$ of $S^T$ takes independently two states, open or
closed, and the probability that $e$ is open is $p \in [0,1]$. The
distribution of states of all edges is denoted by $P_p$. This can be
regarded as a restriction of Bernoulli probability measure $P_p$ on
states of all edges of ${\mathbb Z}^2$ to those of edges of $S^T$.


Also let $S^{T*}$ be the dual graph of $S^T$. Namely, a vertex of
$S^{T*}$ is the central point of a face of $S^T$, and an edge of
$S^{T*}$ is a pair of vertices of $S^{T*}$ such that the
corresponding faces of $S^T$ have a common edge of $S^T$ in their
boundaries. For every edge $e$ of $S^T$ there exists uniquely an
edge $e^*$ of $S^{T*}$ such that $e$ and $e^*$ cross each other. As
usual we say that $e^*$ is open (closed ) if $e$ is open (closed).
Put
\begin{eqnarray*}
\lefteqn{p_c({\mathcal G})}\\
&&= \inf \left\{ p\in [0,1] \ ; \ P_p\bigl[
      \hbox{there exiets an infinite open cluster in ${\mathcal G}$}
                                     \bigr] >0 \right\}
\end{eqnarray*}
for ${\mathcal G}= S^T, G^T$ and $S^{T*}$.
First, we show the following theorem by a standard percolation
argument in section 2.

\begin{theorem}\label{theorem1}
\begin{description}
\item{}{\rm (1)}\quad If $p>p_c(S^T)$, then the infinite cluster is unique a.s.

\item{}{\rm (2)}\quad If $p>p_c(G^T)$, then
$$
\inf_{x,y\in S^T} \tau_p(x,y)>0,
$$
where,
$$
\tau_p(x,y)
  :=\left\{ \hbox{ $x$ and $y$ are in the same open cluster in $S^T$ } \right\} .
$$
\item{}{\rm (3)} \quad If $p<1-p_c(S^{T*})$, then $\tau_p(x,y)$
decays exponentially as $| x-y|\rightarrow \infty $.
\end{description}
\end{theorem}
\begin{remark} {\rm
Kumagai \cite{Kumagai} considered percolation problem on
Sierpi\'nski carpet lattice in a general setting. Namely, let
$T\subset \{0, 1,2, \ldots , L-1\}^2 $, and for $(i,j)\in T$, put
$$
\psi_{i,j}(x_1,x_2)=L^{-1}(i,j)+ L^{-1}(x_1,x_2) \qquad (x_1,x_2)\in
[0,1]^2.
$$
Then there exists an unique closed set $ E^T$ such that
$$ E^T = \bigcup_{(i,j)\in T} \psi_{i,j}(E^T).$$
Assume that $T$ satisfies the following conditions:
\begin{enumerate}
\item $E^T$ is connected,
\item If $(i,j)\in T$ then both $(j,i)$ and $(i,L-1-j)$ are in $T$.
\item $\{ (0,j); 0\leq j\leq L-1\} \subset T$.
\end{enumerate}
Then he proved the uniqueness of the critical probability for the
general Sierpi\'nski carpet lattice $G^T$ generated by $T$, under a
condition related to the crossing probabilities, which we introduce
in section 2. However, it is not clear whether our (original)
Sierpi\'nski carpet lattice $G^T$ or $S^T$ satisfies his condition
and the problem remained open for $T=\{ 0,1,2\}^2\setminus \{(1,1)\}
$.}
\end{remark}

Theorem \ref{theorem1} summarizes where the essential problem lies.
Also the proof of this theorem is a preparation for the argument to
obtain the following

final result.
\begin{theorem}\label{theorem2}
\begin{description}
\item{}{\rm (1)}\quad $p_c(S^T) = p_c(G^T)=1-p_c(S^{T*}).$
\item{}{\rm (2)} \quad The percolation probability
$$
\theta (p) = P_p\left[ \hbox{ there is an open path from the origin
to infinity in
 $S^T$ }\right]
$$
is continuous at $p_c(S^T)$.
\end{description}
\end{theorem}
Combining this theorem with the result in \cite{Sugimine-Takei}, we
obtain that the central limit theorem for the number of open
clusters in $G_n^T$ holds for every $p\in (0,1)$.
\section{Sponge percolation probabilities}

For the proof of our theorems, the sponge percolation probabilities
play important
 roles as usual.
For integers $\ell ,k \geq 1$ and $n\geq 1$, let $G_n(\ell ,k)$ be
the union of shifts of $\ell \times k $ $G_n^T$'s, with $k$ rows
 and $\ell $ columns.
The origin is located at the lower left corner of $G_n(\ell ,k)$. To
be more precise,
$$
G_n(k,\ell ) = \bigcup_{0\leq i \leq \ell -1}\bigcup_{0\leq j\leq
k-1}
 \bigl[ G_n^T + (i\cdot 3^n, j\cdot 3^n) \bigr] .
$$
These may not be subgraphs of $S^T$, but we can consider them as
subgraphs of ${\mathbb Z}^2$, and therefore we can also consider
$P_p$ probabilities of events on edges of these graphs. Let
$A_n(\ell ,k)$ denote the event that there is an open left-right
crossing in $G_n(\ell,
 k)$, and let $B_n(\ell ,k)$ denote the event that there exists an open
up-down crossing
 in $G_n(\ell , k)$.
Let $G_n^*(\ell ,k)$ denote the dual graph of $G_n(\ell ,k) $. By
this, we mean the following graph. 
 First we put a vertex at the center of every finite face
of $G_n(\ell ,k)$, and connect each pair $u^*,v^*$ of these vertices
if the corresponding faces of $G_n(\ell
 ,k)$ have an edge of $G_n(\ell ,k)$ in common on their boundaries. Next,
we add edges of ${\mathbb Z}^2+(\frac{1}{2}, \frac{1}{2})$ if they
connect a finite face of $G_n(\ell ,k)$ to the unique infinite face
of $G_n(\ell ,k)$. Thus, the total edges obtained above form the
edge set of $G_n^*(\ell ,k)$, and the set of all points incident to
some of these edges is the vertex set of $G_n^*(\ell,
 k)$.
We write $A_n^*(\ell ,k)$ for the event that there is a dual closed
left-right crossing
 in
$G_n^*(\ell ,k)$, and $B_n^*(\ell ,k)$ for the event that there is a
dual closed up-down
 crossing in $G_n^*(\ell, k)$.

Let ${\mathbb E}^2$ be the edge set of ${\mathbb Z}^2$. The space of
total edge configurations in ${\mathbb E}^2$ is denoted by
$$\Omega_{{\mathbb E}^2}={ \{ 0,1\} }^{{\mathbb E}^2}.$$
Here, we can define shifts;
$$ \tau_x \omega (b):= \omega (b-x),\qquad x \in {\mathbb Z}^2.$$
Although shifts are not possible to define in the edge configuration
space on $S^T $, $G^T$ or $S^{T*}$, we can define shifts of local
edge events there. Namely, let $A$ be an edge event on a finite
subgraph $V$ of ${\mathbb Z}^2$. If for some $x\in {\mathbb Z}^2$,
$V+x$ is a subgraph of $S^T$, then we can regard $\tau_xA$ as an
edge event in $S^T$. We will use this convention hereafter.

The Kumagai's condition mentioned in the previous section is as
follows.
$$
\limsup_{n\rightarrow \infty }P_p[ A_n(3L,1)]<1 \Leftrightarrow
\limsup_{n\rightarrow \infty }P_p[ A_n(3L,2)]<1
$$
for all $p\in [0,1]$, where in our case, $L=3$. We are not going to
check this condition directly, but investigate crossing
probabilities in a more precise manner.

First, we give an RSW type result which is valid for both original
and dual connection.

\begin{lemma}\label{lemma2.1}
Let $0<a,b<1$. For every $k \geq 1$, there exist continuous
increasing functions $f_k(x), g_k(x,y)$ such that
\begin{description}
\item{}{\rm (1)}\quad If $P_p[ A_n^{(*)}(2,2)]\geq a$, then
$$
P_p\left[ A_n^{(*)}(k,2) \right] \geq f_k(a) .
$$
\item{}{\rm (2)}\quad If $P_p[A_n^{(*)}(2,2)]\geq a$, and
$P_p[A_{n+1}^{(*)}(1,1)] \geq b$, then
$$
P_p\left[ A_{n+1}^{(*)}( k, 1)\right]  \geq g_k(a,b).
$$
\item{}{\rm (3)} \quad $f_k(x)>0$ if $x>0$, and $g_k(x,y)>0$ if $x>0, y>0$.
Further, $f_k(0)=g_k(0,y)=g_k(x,0)=0$ for $(x,y) \in [0,1]^2$, and $
f(1)= g(1,1) =1$.

\end{description}
\end{lemma}

The proof of above lemma is about the same as that of the RSW
theorem  in ${\mathbb E}^2$ (see e.g.\cite{russo}). The functions
$f_k(x)$ can be taken as
\begin{eqnarray*}
f_1(x)&=&f_2(x)=x, \, f_3(x)= (1-\sqrt{1-x})^3, \\
f_{2k+1}(x)&=& x\cdot f_{k+1}(x)\cdot f_{k+2}(x), \\
f_{2k}(x)&=& x \cdot f_{k+1}(x)^2
\end{eqnarray*}
for $k\geq 2 $. As for the functions $g_k$, we give brief sketch of
the proof of (2). We look at $G_{n+1}^{(*)}(2,1)$. By our assumption
and the square root trick, the probability that there is an open (
closed dual ) left-right crossing in $G_{n+1}^T$ passing below $(
2\cdot 3^n, \frac{1}{2}3^{n+1})$ is not less than $ 1- \sqrt{1-b}$.
Let $c$ be the lowest of such a left-right crossing and let $c'$ be
its reflection with respect to the line $\{ x_1=3^{n+1}\} $. Then
the probability that there is an open (closed dual) up-down crossing
in the region above $c\cup c'$ and in $H=G_n(2,3)+(2\cdot 3^n,0)$,
that connects the top side of $H$ with $c$, is not less than
$1-\sqrt{1-f_3(a)}$. Further the probability that there is an open
(closed dual ) left-right crossing in $H$ starting above $( 2\cdot
3^n, \frac{1}{2}3^{n+1})$ is not less than $1-\sqrt{1-a}$. Thus, by
the FKG inequality and by the usual RSW argument the probability
that there is an open ( closed dual ) left-right crossing in
$G_{n+1}^T\cup H$ is not less than
$$
(1-\sqrt{1-b})(1-\sqrt{1-f_3(a)})(1-\sqrt{1-a}).
$$
The same is true for the reflected region of $G_{n+1}^T\cup H$ with
respect to the line $\{ x_1=3^{n+1}\} $. The probability of the
intersection of these two events and the event that there is an open
( closed ) up-down crossing in $H$ can be therefore not less than
$$
(1-\sqrt{1-b})^2(1-\sqrt{1-f_3(a)})^2(1-\sqrt{1-a})^2 \cdot f_3(a),
$$
which we take as $g_2(a,b)$. For general $k\geq 2$, we can take
$$
g_{k+1}(a,b)=b \cdot g_k(a,b)g_2(a,b).
$$

On the other hand, the following lemma is specially for $S^T$ or
$G^T$.

\begin{lemma}\label{lemma2.1.5}
\begin{description}
\item{}{\rm (1)} If $\lim_{n\rightarrow \infty }P_p[ A_n(1,1)]=1$,
then $\lim_{n\rightarrow \infty }P_p[ A_n(k,1)]=1$ for $k\geq 1$.
\item{}{\rm (2)} If $\lim_{n\rightarrow \infty }P_p[ A_n^*(2,2)]=1$, then
$ \lim_{n\rightarrow \infty }P_p[ A_n^*(k,1)]=1$ for $k\geq 1$.
\end{description}
\end{lemma}

\noindent {\it Proof}.\quad (1) \, If there is an open left-right
crossing in $G_n^T$, then this open crossing avoids the central hole
$[3^{n-1},2\cdot 3^{n-1}]^2$. By symmetry and the square root trick,
the probability that there is an open left-right crossing in $G_n^T$
passing below the central hole is not less than
$1-\sqrt{1-P_p[A_n(1,1)]}$. Taking intersection of rotations of this
event and by the FKG inequality, the probability that there is an
open circuit in $G_n^T$ surrounding the central hole
 is not less than
$$ {\left( 1-\sqrt{1-P_p[A_n(1,1)]} \right) }^4.$$
So, with this probability, we can find an open path in
$G_{n-1}(3,1)$ connecting $[0,3^{n-1}]\times \{ 3^{n-1}\} $ with
$[2\cdot 3^{n-1},3^n]\times \{ 3^{n-1}\}$ without touching
$(3^{n-1},2\cdot 3^{n-1})\times \{ 3^{n-1}\} $. Let $E_n$ denote
this event. Then there exists an open left-right crossing in
$G_{n-1}(k,1)$ if
$$ \tau_{(-3^{n-1},0)}E_n \cap E_n \cap \tau_{(3^{n-1},0)}E_n
 \cap \ldots \cap \tau_{((k-2)3^{n-1},0)}E_n $$
occurs. The probability of this event is not less than
$$
P_p[ A_n(k,1)] \geq  {\left( 1-\sqrt{1-P_p[A_n(1,1)]} \right)
}^{4k},
$$
which converges to $1$ if $\lim_{n\rightarrow \infty }P_p[
A_n(1,1)]=1$.

\noindent (2) \, By (1) of Lemma \ref{lemma2.1}, we have
$$
\lim_{n\rightarrow \infty }P_p[ A_n^*(k,2)]=1
$$
for every $k\geq 1$. But for any $k\geq 1$, we always have
$$P_p[A^*_n(3k,2)]\leq P_p(A^*_{n+1}(k,1))$$ by comparison.
Thus we get (2).

\medskip

The crossing probabilities are related to the critical probabilities
introduced in section 1 in the following manner.

\begin{lemma}\label{lemma2.2}
\begin{description}
\item{}{\rm (1)}\quad If $p>p_c(S^T)$, then
$$\lim_{n\rightarrow \infty } P_p[ A_n(2,2) ]=1.$$
\item{}{\rm (2)}\quad If $p>p_c(G^T)$, then
$$\lim_{n\rightarrow \infty }P_p[ A_n(1,1) ]=1,$$
moreover, there exist constants $C>0$ and  $0<\alpha < 1$ such that
$$ P_p\left[ A_n(3,1) \right] \geq 1 - C \alpha^{2^n} $$
for sufficiently large $n$'s.
\item{}{\rm (3)}\quad If $p<1-p_c(S^{T*})$, then
$$\lim_{n\rightarrow \infty }P_p[ A_n^*(2,2)] =1.$$
Further, there exist constants $C^*>0$ and $0<\alpha^*<1$ such that
\begin{equation}\label{eq2.1}
P_p[ A_n^*(3,1)] \geq 1 - C^* (\alpha^*)^{3^n}.
\end{equation}
\end{description}
\end{lemma}

{\it Proof}.\quad (1) \, As in the proof of (2) of the previous
lemma, the probability that there is a closed dual path in
$G_{n-2}(6,2)+(-3^{n-1},2\cdot 3^{n-2})$ connecting the holes
$[-2\cdot 3^{n-1}, 2\cdot 3^{n-2}]\times [3^{n-1},2\cdot 3^{n-1}]$
and $[3^{n-1},2\cdot 3^{n-1}]^2 $ is not less than
$$ {\biggl( 1 - \sqrt{1-f_6(P_p(A_{n-2}^*(2,2)))}\biggr) }^2.$$
Taking intersection of rotations of this event, we see that the
probability that there is a dual closed circuit in $S^T\cap
[-3^n,3^n]^2$ surrounding $[-3^{n-2},3^{n-2}]^2$ is not less than
$$ {\biggl( 1-\sqrt{1-f_6(P_p(A_{n-2}^*(2,2)))}\biggr) }^8. $$
Therefore, if
$$ \limsup_{n\rightarrow \infty }P_p\bigl[ A_n^*(2,2)\bigr]
  =\limsup_{n\rightarrow \infty }( 1- P_p\bigl[ A_n(2,2)\bigr]) >0, $$
then by the second Borel-Cantelli's lemma we find $P_p$-a.s.
infinitely many disjoint closed dual circuits surrounding the
origin. Thus, there is no infinite open cluster in $S^T$, which
implies that $p\leq p_c(S^T)$.

\medskip
\noindent (2)\, If $p>p_c(G^T)$, $P_p$-a.s. for sufficiently large
$n\geq 1$, the infinite open cluster in $G^T$ intersects $G_n^T$.
This open cluster will go out of $G_{n+1}^T$, avoiding the central
hole $[3^n, 2\cdot 3^n]^2$. Therefore it must contain an open
left-right crossing of $G_n^T+(3^n,0)$ or an open up-down crossing
of $G_n^T+(0,3^n)$. By symmetry and the square root trick these
probabilities go to $1$ as $n\rightarrow \infty $. Thus we have
$$ \lim_{n\rightarrow \infty }P_p\bigl[ A_n(1,1)\bigr] = 1. $$
The remaining estimate comes from a simple scaling argument
introduced in \cite{ACCFR}. By symmetry and  the FKG inequality we
have
$$
P_p\left[ A_{n}(9,1)\right] \geq {P_p\left[ A_n(3,1)\right] }^5,
$$
because $A_n(9,1)$ occurs when there is an open
left-right crossing in each $A_n(3,1)+(2j\cdot 3^n,0)$ for
$j=0,1,2,3$ and also there is an open up-down crossing in
$A_n(1,1)+(2j\cdot 3^n,0)$ for $j=1,2,3$. By independence, we have
$$
P_p\left[ A_{n+1}(3,1)\right] \geq \varphi \left( P_p\bigl[
                   A_n(9,1) \bigr] \right) ,
$$
where $\varphi (x) = 1-(1-x)^2$. These together imply that
$$
P_p\left[ A_{n+1}(3,1)\right]
  \geq \psi \left( P_p\bigl[ A_n(3,1)\bigr] \right) ,
$$
where $\psi (x) = 1-(1-x^5)^2$. If $1>x>1-5^{-2}\theta $ for some
$0<\theta <1$, then $\psi (x) \geq 1- 5^{-2}\theta^2$, therefore
once we have
\begin{equation}\label{eq:scaling}
P_p\bigl[ A_n(3,1)\bigr] \geq 1-5^{-2}\theta ,
\end{equation}
for some $n$, we have
$$ P_p\bigl[ A_{n+k}(3,1) \bigr] \geq 1-5^{-2}\theta^{2^k} $$
for all $k\geq 0$ and the statement of (2) is true since by Lemma
\ref{lemma2.1.5}, (1),
$$
P_p\left[ A_n(3,1)\right]\rightarrow 1.
$$

\medskip
\noindent The proof of (3) is done by usual argument( see e.g,
Theorem 8.97 of \cite{grimmett}). Namely, if $p<1-p_c(S^{T*})$, then
$P_p$-a.s. there is an infinite closed dual cluster. Putting
$$
\theta^*(p)=P_p\left[
          \hbox{one of $(\pm\frac{1}{2},\pm\frac{1}{2})$ is in an }
          \hbox{infinite closed dual cluster}
 \right] ,
$$
we have by symmetry and the FKG inequality
$$
P_p\biggl[ A_n^*(2,2)\biggr] \geq
  {\biggl( \frac{1}{4}\theta^*(p)\biggr) }^2(1-p)^4.
$$
Then by Lemma\ref{lemma2.1} and by the same argument in the proof of
(1), we have $P_p$-a.s. for every finite $\Lambda \ni 0$, there is a
closed dual circuit surrounding $\Lambda $. Since an infinite closed
dual cluster intersects large $\Lambda \ni 0$ with probability close
to $1$, by symmetry and the FKG inequality we have
$$
P_p\biggl[ A_n^*(2,2)\biggr] \rightarrow 1.
$$
As for the exponential estimate (\ref{eq2.1}), note that
$A_{n+1}^*(3,1)$ occurs if there exists a closed left-right crossing
 in one of  $G_n(9,1), G_n(9,1)+(0,3^n)$ and $G_n(9,1)+(0,2\cdot 3^n)$.
Each probability is not less than $P_p[A_n^*(9,1)]$, and hence we
have desired estimate by the same scaling argument as in the proof
of (2) of the lemma.

\bigskip

\noindent {\it Proof of Theorem\ref{theorem1}}.\quad (1) Let
$p>p_c(S^T)$. Then by Lemma\ref{lemma2.2},
$$
\lim_{n\rightarrow \infty }P_p[ A_n(2,2)]=1.
$$
If there are more than one infinite open clusters $P_p$-a.s., then
there exists an infinite closed dual cluster separating different
infinite open clusters $P_p$-a.s., which implies that
$$
\lim_{n\rightarrow \infty }P_p[ A_n^*(2,2)]=1
$$
by the same argument as in the proof of Lemma\ref{lemma2.2}, (3).
This is a contradiction, since $P_p[A_n^*(2,2)]=1-P_p[A_n(2,2)]$.

\noindent (2) It is sufficient to prove that the inequality
$p>p_c(G^T)$ guarantees that
\begin{equation}\label{eq,tau}
\inf_{x \in G^T} P_p\bigl[
 \hbox{ $x$ is connected to the origin by an open path in $G^T$}
   \bigr] >0.
\end{equation}
Let $p>p_c(G^T)$. Then by Lemmas \ref{lemma2.1.5} and
\ref{lemma2.2}, we have
$$
\lim_{n\rightarrow \infty }P_p[ A_n(3,1)]=1.
$$
Put
$$\Delta_n=A_{n-1}(3,1)\cap B_{n-1}(1,3)\cap
  \tau_{(0,2\cdot 3^{n-1})}A_{n-1}(3,1)\cap
   \tau_{(2\cdot 3^{n-1},0)}B_{n-1}(1,3).
$$
Then by the FKG inequality we have
$$
P_p[\Delta_n]\geq {\left( P_p[ A_{n-1}(3,1)]\right) }^4.
$$
On $\Delta_n$, there is a big open cluster in $G_n^T$ such that
every open crossing in $G_{n-1}(3,1), G_{n-1}(1,3),
G_{n-1}(3,1)+(0,2\cdot 3^{n-1})$ and $G_{n-1}(1,3)+(2\cdot
3^{n-1},0)$ connecting shorter sides of the rectangle is a subset of
this big open cluster. We call this open cluster as the spanning
cluster in $G_n^T$. Take a subgraph of $G_n^T$ such that it is a
shift of $G_{n-1}^T$. Then it is written as
$$
G_{n-1}^T+3^{n-1}(i,j)
$$
 for some $(i,j)\in T$.
Note that the spanning cluster of $G_n^T$ includes the spanning
cluster of $G_{n-1}^T+3^{n-1}(i,j)$ if the event
$\tau_{(i3^{n-1},j3^{n-1})}\Delta_{n-1}\cap \Delta_n$ occurs. For
$x\in G^T$, we define the level of $x$ as $ \ell (x) = \inf\{ n : x
\in G_n^T\} $. Also we define $G_n^T(x)$, the n-th box of $x$, by
the shift of $G_n^T$ such that it is a subgraph of $G^T$ and it
contains $x$. If there are more than one such boxes, then we take
the nearest one to the origin as $G_n^T(x)$. For convenience, let us
write
$$ G_n^T(x) = G_n^T + w_n(x).$$
Note that $w_n(0)=0$ for every $n\geq 1$. Since $p>p_c(G^T)$, by (2)
of Lemma\ref{lemma2.2}, we have
$$
P_p[ A_n(3,1)]\geq 1 - C\alpha^{2^n}>0\qquad n\geq m_0
$$
if $m_0$ is sufficiently large. If $\ell (x)\leq m_0$, then we have
$$
\tau_p(0,x)\geq P_p[ \mbox{ all edges of $G_{m_0}^T$ are open
}]:=\xi (p)>0.
$$
If $\ell (x)>m_0$, then consider the event
$$
\bigcap_{n=m_0+1}^{\ell (x)} \Delta_n \cap \bigcap_{n=m_0+1}^{\ell
(x)}\tau_{w_n(x)}\Delta_n \cap \{ \mbox{ all edges of $G_{m_0}^T\cup
G_{m_0}^T(x)$ are open }\} .
$$
On this event $x$ and the origin are connected by an open path in
$G^T$, and by the FKG inequality this probability is bounded by
$$
{\left( \xi (p) \prod_{n=m_0+1}^\infty (1-C\alpha^{2^n}) \right)
}^2>0
$$
from below.

\noindent (3) Let $p<1-p_c(S^{T*})$ and let $x,y \in S^T$. We define
$n(x,y)$ by
$$
n(x,y) = \left\{
   \begin{array}{ll}
    \max\{ n: d(G_n^T(x), G_n^T(y))>0 \} ,  & \hbox{ if }
            |x-y|>6\sqrt{2}, \\
    0, & \hbox{ if } |x-y|\leq 6\sqrt{2}.
   \end{array}
         \right.
$$
Then we have
$$|x-y|\leq 6\sqrt{2}\cdot 3^{n(x,y)}$$
Now consider the following 4 rectangles surrounding
$G_{n(x,y)}^T(x)$.
\begin{eqnarray*}
R_n^1(x)&=& [-3^n, 2\cdot 3^n]\times [3^n, 2\cdot 3^n]+w_n(x),\\
R_n^2(x)&=& [-3^n, 2\cdot 3^n]\times [-3^n, 0]+w_n(x),\\
R_n^3(x)&=& [3^n, 2\cdot 3^n]\times [-3^n, 2\cdot 3^n]+w_n(x),\\
R_n^4(x)&=& [-3^n, 0]\times [-3^n, 2\cdot 3^n]+w_n(x),
\end{eqnarray*}
where $n=n(x,y)$. Let $\Delta_{n(x,y)}^*(x)$ be the event that there
exists a closed dual crossing in each of above rectangles,
connecting shorter sides of each $R_n^i(x), \ i=1,2,3,4$. Then the
probability of  $\Delta_{n(x,y)}^*(x)$ is not less than
$$
{\biggl( 1 - C^*\cdot (\alpha^*)^{3^n}\biggr) }^4,
$$
by Lemma\ref{lemma2.2}. Finally, note that there is a closed dual
circuit surrounding $x$, such that $y$ is outside of this circuit.

\bigskip
Now we proceed to the proof of Theorem\ref{theorem2}. Let
$C_{n,b\ell }$ be the event that satisfies all of the following
conditions;

\begin{enumerate}
\item In the rectangle
$$ G_n (6, 2)+ ( - 3^{n+1}, 2\cdot 3^n ),$$
there exists an open left-right crossing which ends at the boundary
of the central hole $ [3^{n+1}, 2\cdot 3^{n+1}]^2 $ of $G^T_{n+2}$,
\item in the rectangle
$$ G_n (2,3)+(2\cdot 3^n,0), $$
there exists an open up-down crossing which ends at the boundary of
the central hole of $G^T_{n+2}$, and
\item any pair of above open crossings are connected by an open path
in
$$ G^T_{n+2}\cap \biggl\{
    [G_n(6,2)+ (-3^{n+1}, 2\cdot 3^n )] \cup
   [G_n(2,3)+ (2\cdot 3^n,0)] \biggr\} .
$$
The symbol ``$b\ell$ `` stands for ``bottom and left''.
\end{enumerate}

The key to the proof of the first equality of (1) of Theorem
\ref{theorem2} is the following lemma.

\begin{lemma}\label{lemma2.3}
If $\lim_{n\rightarrow \infty } P_p[ A_n(2,2)]=1$, then
$$\lim_{n\rightarrow \infty }P_p[ C_{n,b\ell }]=1. $$
\end{lemma}

We will give proof of this lemma in the next section. By this lemma,
we can show that

\begin{lemma}\label{lemma2.4}
If $\lim_{n\rightarrow \infty }P_p[ A_n(2,2)]=1$, then
$$ \lim_{n\rightarrow \infty }P_p[ A_n(1,1)] = 1, $$
which means that $p\geq p_c(G^T)$.
\end{lemma}

Combining this lemma with  (1) of Lemma \ref{lemma2.2}, we obtain
the equality
$$p_c(S^T)=p_c(G^T).$$

\bigskip
\noindent {\it Proof of Lemma \ref{lemma2.4}}.\quad Let $C_{n,br}$
be the reflected event of $C_{n,b\ell }$ with respect to the line
$\{ x^1 = \frac{ 3^{n+2}}{2} \} $, and $C_{n,t\ell }$ be the
reflected event of $C_{n, b\ell }$ with respect to $\{
x^2=\frac{3^{n+2}}{2} \} $. Also, we define $C_{n,tr}$ as the
reflected event of $C_{n,b\ell }$
 with respect to the point $(\frac{3^{n+2}}{2},\frac{3^{n+2}}{2})$.
Now, let
$$C_{n, all} := C_{n,b\ell }\cap C_{n,br}\cap C_{n,t\ell }\cap C_{n,tr}$$
and let
$$C'_{n,all}:= C_{n,all}+(3^{n+2},0).$$
Let $D_n$ be the event that there is an open left right crossing in
$G^T_{n}$. Then $\omega \in D_{n+2}\cap C_{n,all}$ implies that
there is an open left right crossing in the rectangle
$$ [-3^{n+1},3^{n+2}+3^{n+1}]\times [0, 3^{n+2}]. $$
Let
$$
C_{n,bottom}= C_{n,b\ell }\cap C_{n,br}\quad C_{n,top}=C_{n,t\ell
}\cap C_{n,tr}.
$$
Then we have
\begin{equation}\label{eq,7}
\{ C_{n,top}\cap C''_{n-1,all}\cap D''_{n+1} \} \cup \{
C_{n,bottom}\cap C'_{n-1,all}\cap D'_{n+1} \} \subset D_{n+2},
\end{equation}
where $C''_{n-1,all}$ and $D''_{n+1}$ are shifts of $C_{n-1,all}$
and $D_{n+1}$ by $(3^{n+1},2\cdot 3^{n+1})$, and $C'_{n-1,all}$ and
$D'_{n+1}$ are shifts of $C_{n-1,all}$ and $D_{n+1}$ by
$(3^{n+1},0)$. Since $C_{n,top}\cap C''_{n-1,all}\cap D''_{n+1}$ and
$C_{n,bottom}\cap C'_{n-1,all}\cap D'_{n+1}$ are independent,
(\ref{eq,7}) implies the following inequality
\begin{equation}\label{eq,8}
P_p[D_{n+2}]\geq f(P_p[ C_{n,bottom}\cap C_{n-1,all}\cap D'_{n+1}]),
\end{equation}
where $f(t)=2t-t^2$, which is an increasing function in the interval
$[0,1]$.

By the FKG inequality and by the equality $ P_p[D'_{n+1} ] =
P_p[D_{n+1}]$, we have
\begin{equation}\label{eq,9}
P_p[ C_{n,bottom}\cap C_{n-1,all}\cap D'_{n+1}] \geq
P_p[C_{n,bottom}]P_p[C_{n-1,all}]P_p[D_{n+1}]
\end{equation}

By Lemma \ref{lemma2.3}, If $\lim_{n\rightarrow \infty} P_p[
A_n(2,2)]=1$, then we have
$$
\lim_{n\rightarrow \infty} P_p[C_{n,bottom}]
 =\lim_{n\rightarrow \infty } P_p[ C_{n,all}]=1.
$$
Therefore for every $\varepsilon >0$, there exists some $N>0$ such
that
\begin{equation}\label{eq,10}
P_p[C_{n,bottom}]\geq P_p[C_{n,all}]> 1-\varepsilon
\end{equation}
for every $n>N$. Combining (\ref{eq,8}--\ref{eq,10}), we have
\begin{equation}\label{eq,11}
P_p[D_{n+2}]\geq f((1-\varepsilon)^2P_p[D_{n+1}])
\end{equation}
for every $n>N+1$. This implies that
$$
\liminf_{n\rightarrow \infty }P_p[D_n] \geq x_{\varepsilon },
$$
where $x_{\varepsilon }$ is the unique solution to
$$
t = f( (1-\varepsilon )^2 t ),
$$
which converges to $1$ as $\varepsilon \rightarrow 0$. Since $D_n =
A_n(1,1)$, we have
$$
\lim_{n\rightarrow \infty }P_p[A_n(1,1)]=\lim_{n\rightarrow \infty
}P_p(D_n)=1.
$$
By Lemma\ref{lemma2.1.5} and by the proof of Theorem
\ref{theorem1},(2), this implies the inequality (\ref{eq,tau}). But
then we have
\begin{eqnarray*}
\lefteqn{P_p[\hbox{ the open cluster of the origin is an infinite cluster}]}\\
&&=\lim_{n\rightarrow\infty}
   P_p\left[
  \begin{array}{l}
   \hbox{ the origin is connected to $\{ x_1=3^n \mbox{ or }x_2=3^n\} $}\\
    \hbox{ by an open path in $G_n^T$}
  \end{array}
   \right] \\
&&\geq \liminf_{n\rightarrow \infty}P_p[
     \hbox{ the origin is connected to $(3^n,0)$ by an open path in $G_n^T$}
      ] .
\end{eqnarray*}
The right hand side of the above inequality is positive by
(\ref{eq,tau}).

\bigskip

The proof of the second equality in (1) of Theorem\ref{theorem2} and
the proof of (2) of Theorem\ref{theorem2} are postponed to the last
section.

\section{ Branching argument }

In this section, we prove Lemma\ref{lemma2.3}. Before going into the
detail, we give rough idea of the proof. By the condition that
$$\lim_{n\rightarrow \infty }P_p[ A_n(2,2)]=1$$
together with Lemma\ref{lemma2.1}, with high probability there is an
open up-down crossing $\gamma $ in $G_n(2,3)+(2\cdot 3^n, 0)$.
Further, with high probability we can find an open path branching
out from this open up-down crossing in $G_{n-1}(6,2)+(2\cdot 3^n,
5\cdot 3^{n-1})$ to the line
 $\{ x_1=2\cdot 3^n\} $.
In the same way, with high probability we can find an open path
branching out from the open left-right crossing $\delta $ of
$G_n(3,2)+( 0, 2\cdot 3^n)$ to the line $\{ x_2=2\cdot 3^n\} $
 in $G_{n-1}(2,6)+(5\cdot 3^{n-1}, 2\cdot 3^n)$.
From these branches we can find open paths branching out with high
probability. By these branches, the possibility of connecting the
original up-down crossing $\gamma
 $
and the left-right crossing $\delta $ increases. As we keep on this
procedure, we can find with high probability many branches of
original open paths, which become closer

and closer. Therefore with high probability, we can find many pairs
of open branches of $\gamma $ and $\delta $ which are very close.
Finally, connecting one of such pair of branches costs loss of only
small probability.

In the actual procedure, we have to choose $\gamma $ and $\delta $
so that they are also close to each other. For this, we will use
site percolation on a rooted binary tree. Now, let us begin with
some notations. By $\mathbb{T}_2$, we mean a rooted binary tree. The
origin of $\mathbb{T}_2$ is denoted by $\mathbf{0}$. A point of
$\mathbb{T}_2\setminus \{ \mathbf{0}\} $
 is denoted by $\mathbf{j}=(j_1, \ldots , j_n)$
with $j_1, \ldots , j_n \in \{ 1,2\} $. The point $(1)$ is the first
child of $\mathbf{0}$, and $(1,2)$ is the second child of $(1)$, and
so on.

Let $N$ be sufficiently large and fixed. We will  specify later how
large $N$ should be. Let
\begin{equation}\label{eq,12}
V_0^{(N)}=[-3^{N}, 3^{N}]\times [-3^N, 2\cdot 3^{N+1}+ 3^{N-1}]
\end{equation}
and
\begin{equation}\label{eq,13}
J_0^{(N)}=[-3^N-3^{N-1}, 3^N +3^{N-1}]\times
  [2\cdot 3^{N+1}-3^{N-1}, 2\cdot 3^{N+1}+ 3^{N-1}]
\end{equation}
and set
\begin{equation}\label{eq,14}
{\tilde V}_0^{(N)} =
      V_0^{(N)}\cup  J_0^{(N)}.
 \end{equation}
This is the mother shape for the straight connection. Further, we
introduce
\begin{equation}\label{eq,15}
\begin{array}{lcl}
J_{0,\ell }^{(N)}&=&[-3^N-3^{N-1},-3^N+3^{N-1}]\times
    [2\cdot 3^{N+1}-3^{N-1}, 2\cdot 3^{N+1}+3^{N-1}], \\
J_{0,r}^{(N)}&=&[3^N-3^{N-1}, 3^N+3^{N-1}]\times
    [2\cdot 3^{N+1}-3^{N-1}, 2\cdot 3^{N+1}+3^{N-1}].
\end{array}
\end{equation}
For the branching connection, the mother shape is different. Let
\begin{equation}\label{eq,16}
\Lambda_0^{(N)}= [-3^{N+1}-3^{N-1}, 3^N]\times [-3^N, 3^N],
\end{equation}
\begin{equation}\label{eq,17}
I_0^{(N)}= [-3^{N+1}-3^{N-1},-3^{N+1}+3^{N-1}]\times [ -3^N-3^{N-1},
3^N+3^{N-1}]
\end{equation}
and
\begin{equation}\label{eq,18}
{\tilde \Lambda }_0^{(N)}= \Lambda_0^{(N)}\cup I_0^{(N)},
\end{equation}
which is the mother shape for the branching connection. Further, we
introduce
\begin{equation}\label{eq,19}
\begin{array}{lcl}
I_{0,t}^{(N)}&=&[ -3^{N+1}-3^{N-1}, -3^{N+1}+3^{N-1}]\times
   [3^N-3^{N-1}, 3^N+3^{N-1}]\\
I_{0,b}^{(N)}&=&[-3^{N+1}-3^{N-1}, -3^{N+1}+3^{N-1}]\times
   [-3^N-3^{N-1}, -3^N+3^{N-1}].
\end{array}
\end{equation}
The scaled shapes ${\tilde V}_k^{(N)}$ and ${\tilde \Lambda
}_k^{(N)}$ are defined
 by
$$
{\tilde V}_k^{(N)}= 3^{-k}{\tilde V}_0^{(N)},\quad {\tilde \Lambda
}_k^{(N)}=3^{-k}{\tilde \Lambda }_0^{(N)}.
$$
In the same way, we define
$$
V_k^{(N)}=3^{-k}V_0^{N)}, \ J_k^{(N)}=3^{-k}J_0^{(N)}, \
J_{k,*}^{(N)}=3^{-k}J_{0,*}^{(N)}\quad *= \ell, r,
$$
and
$$
\Lambda_k^{(N)}=3^{-k}\Lambda_0^{(N)}, \ I_k^{(N)}=3^{-k}I_0^{(N)},
\ I_{k,**}^{(N)}=3^{-k}I_{0,**}^{(N)}, \quad **=t, b.
$$
Let $\theta $ denote the rotation of 90 degrees with respect to the
origin. We put
\begin{eqnarray*}
B_{\mathbf 0}^{(N)}&=& {\tilde V}_0^{(N)}+(3^{N+2},0) \\
B_{\mathbf 0}^{(N)\dagger } &=& \theta^{-1}{\tilde
V}_0^{(N)}+(0,3^{N+2})
\end{eqnarray*}
and ${\mathcal B}^{(N)}_0=B_{\mathbf 0}^{(N)}\cup B_{\mathbf
0}^{(N)\dagger }$. Note that $B_{\mathbf 0}^{(N)\dagger }$ is the
reflection of $B_{\mathbf 0}^{(N)}$ with respect to the line $\{
x_1=x_2\} $, and that ${\mathcal B}^{(N)}_{\mathbf 0}$ is a subgraph
of $S^T$. Let ${\mathbf i}=(i_1, \ldots ,i_n)\in {\mathbb
T}_2\setminus \{ {\mathbf 0}\} $. We introduce the following
notations for ${\mathbf i}$.
\begin{eqnarray*}
| {\mathbf i }| &=& n = \hbox{ the generation that ${\mathbf i}$
belongs },
  \quad |{\mathbf 0}|=0,\\
N_2({\mathbf i})&=& \# \{ \alpha \in \{ 1,\ldots , n\} \ ; i_\alpha
=2 \},
\quad N_2({\mathbf 0})=0 \\
\epsilon ({\mathbf i})&=& N_2({\mathbf i})\quad  ({\rm mod}\  2 )
\end{eqnarray*}
\begin{equation}\label{eq,20}
\left\{
\begin{array}{lll}
\tau_0 &=& 0,  \\
 && \ldots  \\
\tau_\nu &=& \left\{
           \begin{array}{l}
              \min\bigl\{ \tau_{\nu -1}< \alpha \leq n \ ;
              i_\alpha =2 \bigr\} , \\
          \infty \hbox{ if the above set is empty },
        \end{array}
   \right.
\end{array}
\right.
\end{equation}
\begin{equation}\label{eq,21}
 \left\{
\begin{array}{lll}
\ell_v({\mathbf 0})&=& 2\cdot 3^{N+1} \\
\ell_h({\mathbf 0})&=& 0\\
\ell_v({\mathbf i})&=&
    \sum_{\mu = 0}^{\lfloor \frac{N_2({\mathbf i})}{2}\rfloor}\
    \sum_{\alpha =\tau_{2\mu }}^{(\tau_{2\mu +1}-1)\wedge n}
     2 \cdot 3^{N-\alpha +1}, \\
\ell_h({\mathbf i})&=&
   \sum_{\mu = 1}^{\lfloor \frac{N_2({\mathbf i})+1}{2}\rfloor}\
    \sum_{\alpha =\tau_{2\mu -1} }^{(\tau_{2\mu }-1)\wedge n}
     2 \cdot 3^{N-\alpha +1},
\end{array}
\right.
\end{equation}
where $\lfloor x \rfloor $ denotes the largest integer not larger
than $x$, and $a\wedge b = \min \{ a, b\} $.

For a point ${\mathbf j}=(j_1,\ldots ,j_n)\in {\mathbb T}_2\setminus
\{ {\mathbf 0}\}
 $,
 we define vectors $x({\mathbf j}), \ x^\dagger ({\mathbf j})$ and the sets
${\mathcal B}^{(N)}_{\mathbf j}=B^{(N)}_{\mathbf j}\cup
B^{(N))\dagger }_{\mathbf j}$

in the following way.
$$
x({\mathbf j})= \left\{
   \begin{array}{ll}
    \bigl(3^{N+2}- \ell_h({\mathbf j}^*), \ell_v({\mathbf j}^*) \bigr)
          & \hbox{ if $j_n=1$ },\\

    \left(
           3^{N+2}-\ell_h({\mathbf j}^*)- \delta_h({\mathbf j}^*),
                   \ell_v({\mathbf j}^*)+ \delta_v({\mathbf j}^*)
            \right) & \hbox{ if $j_n=2$},
  \end{array}
\right.
$$
where ${\mathbf j}^*$ is the parent of ${\mathbf j}$, i.e. ${\mathbf
j}^*= (j_1, \ldots ,j_{n-1})$, and
\begin{eqnarray*}
\delta_v({\mathbf 0})=0,  \quad
\delta_v({\mathbf j}) &=& 3^{N-|{\mathbf j}|}\epsilon ({\mathbf j}),\\
\delta_h({\mathbf 0})=3^N, \quad \delta_h({\mathbf j}) &=&
      3^{N-|{\mathbf j}|}(1- \epsilon ({\mathbf j})),
\end{eqnarray*}
and we write $x^\dagger ({\mathbf j})$ for the symmetric point of
$x({\mathbf j})$ with respect to the line $\{ x_1=x_2 \} $.

\begin{enumerate}
\item If $j_n =1$, put
$$
B^{(N)}_{\mathbf j}
 =  \theta^{\epsilon ({\mathbf j}^*)}{\tilde V}^{(N)}_n
         + x({\mathbf j}).
$$
\item If $j_n=2$,  put
$$
B^{(N)}_{\mathbf j}
 =  \theta^{-\epsilon ({\mathbf j}^*)}{\tilde \Lambda }^{(N)}_n
          + x({\mathbf j}).
$$
\end{enumerate}
We define $B^{(N)\dagger }_{\mathbf j}$ as the reflection of
$B^{(N)}_{\mathbf j}$ with respect to the line $\{ x_1=x_2\} $.

Note that ${\mathcal B}_{\mathbf j}^{(N)}$ is a subgraph of $S^T$.
For a rectangle $R \subset S^T$, we say that there exists an open
traversing in $R $ if there exists an open path in $R $ which
connects shorter sides of $R $. Let
$$
T(R )=\left\{
       \hbox{ there exists an open traversing in }R
              \right\} .
$$
Further, we define
\begin{eqnarray*}
S(J_n^{(N)})&=& T(J_n^{(N)}) \cap T(J_{n,\ell }^{(N)}) \cap T(J_{n, r}^{(N)} )\\
S({\tilde V}_n^{(N)})&=& T(V_n^{(N)})\cap S(J_n^{(N)}),
\end{eqnarray*}
where the traversing direction of $T(J_{n,*}^{(N)})$ is chosen to be
the same direction as the traversing of $V_n^{(N)}$. Similarly, let
\begin{eqnarray*}
S(I_n^{(N)})&=&T(I_n^{(N)})\cap T(I_{n, t}^{(N)})\cap T(I_{n, b}^{(N)}), \\
S({\tilde \Lambda }_n^{(N)})&=& T(\Lambda_n^{(N)})\cap S(I_n^{(N)}),
\end{eqnarray*}
where the traversing direction of $T(I_{n,*}^{(N)})$ is chosen to be
the same direction
 as
the traversing of $\Lambda_n^{(N)}$. Note that these are all edge
events of $S^T$. Let $\theta $ denote the induced transformation on
${\mathbb E}^2$ by the rotation $\theta $, i.e.,
$$\theta \omega (b)= \omega (\theta^{-1} b) \qquad b \in {\mathbb E}^2.$$
Let
$$
S( B_{\mathbf 0}^{(N)})= \tau_{(3^{N+2},0)}S({\tilde V}_0^{(N)}),
\quad S( B_{\mathbf 0}^{(N)\dagger })=
  \tau_{(0,3^{N+2})}\theta^{-1} S({\tilde V}_0^{(N)}),
$$
and for ${\mathbf j}\in {\mathbb T}_2\setminus \{ \mathbf{ 0}\} $,
we define an edge event $S( B_{\mathbf j}^{(N)})$ by
$$
S( B_{\mathbf j}^{(N)}) =\left\{
     \begin{array}{ll}
     \tau_{x({\mathbf j})}\theta^{\epsilon ({\mathbf j}^*)}S({\tilde V }_{|{\mathbf
 j}|}^{(N)}), & \hbox{ if } {\mathbf j}=( {\mathbf j}^*, 1),\\
     \tau_{x({\mathbf j})}\theta^{-\epsilon ({\mathbf j}^*)}S({\tilde \Lambda}_{|{\mathbf
 j}|}^{(N)}), &\hbox{ if }{\mathbf j}=( {\mathbf j}^*,2),
    \end{array}
\right.
$$
where ${\mathbf j}^*$ is the parent of ${\mathbf j}$. This is
actually an edge event on $ B_{\mathbf j}^{(N)}$. In the same way,
we define an edge event $S(B_{\mathbf j}^{(N)\dagger })$ by the
reflected event of $S( B_{\mathbf j}^{(N)})$ with respect to the
line $\{ x_1=x_2\} $, i.e.,
$$
S(B_{\mathbf j}^{(N)\dagger })=\left\{
  \begin{array}{ll}
    \tau_{x^\dagger ({\mathbf j})}
     \theta^{-1-\epsilon ({\mathbf j}^*)}
     S({\tilde V}_{|{\mathbf j}|}^{(N)}),
        &\hbox{ if } {\mathbf j}=({\mathbf j}^*,1),\\
   \tau_{x^\dagger ({\mathbf j})}
     \theta^{1+\epsilon ({\mathbf j}^*)}
     S({\tilde \Lambda }_{|{\mathbf j}|}^{(N)}),
       &\hbox{ if } {\mathbf j}=({\mathbf j}^*,2).
  \end{array}
\right.
$$
For convenience, we introduce site variables $X(t), X^\dagger (t)$
and $Z(t)$ for $t \in {\mathbb T}_2$ by
$$
X(t)=\left\{
  \begin{array}{ll}
   1, &\hbox{ if } S( B_t^{(N)}) \hbox{ occurs,}\\
   0, &\hbox{ otherwise, }
  \end{array}
   \right.
$$
$$
X^\dagger (t)=\left\{
  \begin{array}{ll}
   1, &\hbox{ if } S( B_t^{(N)\dagger }) \hbox{ occurs,}\\
   0, &\hbox{ otherwise, }
  \end{array}
   \right.
$$
and
$$
Z(t) = X(t)X^\dagger (t).
$$
Let us start with simple facts that can be derived from the
assumption of Lemma \ref{lemma2.3}. By Lemma \ref{lemma2.1} and the
assumption of Lemma \ref{lemma2.3}, we know that for every
$\varepsilon >0$, we can find  $m_0\geq 1$ such that
\begin{equation}\label{eq,22}
P_p\biggl[ A_n(8,2) \biggr] \geq 1-\varepsilon
\end{equation}
for every $n\geq m_0$. Therefore if $ N\geq m\geq m_0+1$, then by
the FKG inequality we have
\begin{eqnarray*}
P_p\biggl[ S({\tilde V}_{N-m}^{(N)})\biggr] &\geq &(1-\varepsilon )^4,\\
P_p\biggl[ S({\tilde \Lambda }_{N-m}^{(N)})\biggr] &\geq &
(1-\varepsilon )^4.
\end{eqnarray*}
This means that for $N\geq m \geq m_0+1$,
$$
P_p\bigl[ Z(t)=1  \bigr] \geq ( 1- \varepsilon )^8
$$
for every $t \in {\mathbb T}_2$ with $|t|\leq N-m$. Note that $\{
X_{\mathbf j}, |{\mathbf j}|\leq N-m \} $ and $\{ X^\dagger_{\mathbf
j}, |{\mathbf j}|\leq N-m \} $ are independent. Further it is easy
to see that
$$
P_p[ Z_{\mathbf j}=1 \mid Z_t = \varepsilon_t, |t-{\mathbf j}|>1 ]
  = P_p[ Z_{\mathbf j}=1 ] \geq (1-\varepsilon )^8,
$$
where $|t-{\mathbf j}|$ is the graph distance of $t$ and ${\mathbf
j}$ in ${\mathbb T}_2$. In this sense $Z_{\mathbf j}$'s are
$1$-dependent. Then by \cite{Liggett} , p.14, Theorem B26, the
distribution of $\{ Z(t), |t|\leq N-m\} $ dominates that of
Bernoulli random variables $\{ W(t), |t|\leq N-m \} $ with
$$ P[ W(t)=1 ]=p_\varepsilon , $$
where $p_\varepsilon $ is given by the unique positive solution to
$$ 1-(1-\sqrt{p})^4 = (1-\varepsilon )^8.$$
Note that $p_\varepsilon \leq (1-\varepsilon )^8$ and that
$p_\varepsilon \rightarrow 1$ as $\varepsilon \rightarrow 0$.

\medskip
\noindent {\it Proof of Lemma \ref{lemma2.3}}

Let us fix an integer $m$ with $ m\geq m_0+1$ and take $N\geq m$.
 By the above observation, we can construct $0$-$1$ valued
random variables $\{ {\tilde Z}(t), W(t); t\in {\mathbb T}_2,
  |T|\leq N-m \} $ on a probability space
 $({\tilde \Omega }, {\tilde F}, {\tilde P})$  such that
\begin{enumerate}
\item $ \{ W(t); t\in {\mathbb T}_2, |t|\leq N-m \} $ is i.i.d. with
$${\tilde P}( W(t) =1 ) = p_\varepsilon , $$
\item the distribution of
$\{ {\tilde Z}(t) ; t\in {\mathbb T}_2, |t|\leq N-m\} $ is the same
as that of $\{ Z(t) ; t\in {\mathbb T}_2, |t|\leq N-m \} $,
\item $\displaystyle {\tilde P}( C_{\mathbf 0}({\tilde Z}) \supset C_{\mathbf 0}(W)
 )=1, $
where $C_{\mathbf 0}({\tilde Z})$ and $C_{\mathbf 0}(W)$ are open
clusters of ${\mathbf 0}$ in the configurations $\{ {\tilde Z}(t) ;
t\in {\mathbb T}_2, |t|\leq
 N-m \} $
and $\{ W(t) ; t\in {\mathbb T}_2, |t|\leq N-m \} $, respectively.
\end{enumerate}

Let
$${\mathcal Z}_n (W)= \# \{ {\mathbf j}\in {\mathbb T}_2\setminus \{ {\mathbf 0}\}
        ; \ |{\mathbf j}|=n , {\mathbf j}\in C_{\mathbf 0}(W)\} . $$

Then conditioned that $W({\mathbf 0})=1$, ${\mathcal Z}_n(W)$ is a
Golton-Watoson branching process with offspring distribution;
$$
p_0 = {\bigl( 1- p_\varepsilon \bigr) }^2,\
  p_1= 2p_\varepsilon \bigl( 1- p_\varepsilon \bigr) , \
  p_2= p_\varepsilon^2, \  p_k=0 \hbox{ for } k\geq 3.
$$
If $0<\varepsilon$ is sufficiently small, then this branching
process is supercritical and for any integer $k\geq 1$,
$$
{\tilde P}( {\mathcal Z}_n(W)\geq k ) \rightarrow 1 - q_\varepsilon
$$
as $n$ goes to infinity, where $q_\varepsilon $ is the extinction
probability of ${\mathcal Z}_n(W)$, which goes to $0$ as
$\varepsilon \rightarrow 0$ ( cf. \cite{Harris}, p.8, Theorem
I.6.1). Let $N_k\geq m$ be so large that
$$
{\tilde P}\bigl( {\mathcal Z}_{N-m}(W)\geq k \bigr) \geq 1 -
2q_\varepsilon
$$
for every $N\geq N_k$. Let
$$
{\mathcal B} = \bigcup\left\{
                      {\mathcal B}_t^{(N)} \, ; \, t \in {\mathbb T}_2,
                    |t|\leq N-m
                       \right\} ,
$$
and fix an integer $N\geq N_k$ and  a configuration on ${\mathcal
B}$ such that
$$
C_{\mathbf 0}(Z)\cap \{ t\in {\mathbb T}_2 ; |t| = N-m \} \not=
\emptyset .
$$
This event occurs in ${\mathcal B}$. We take a point ${\mathbf j}\in
C_{\mathbf 0}(Z)$ with $|{\mathbf j}|=N-m$. Then the unique path
$\xi $ in ${\mathbb T}_2$ which connects ${\mathbf j}$ with
${\mathbf 0}$ is included in $C_{\mathbf 0}(Z)$, therefore there
exists an open path  $\gamma_1$ in $\cup_{t \in \xi }B_t^{(N)}$ that
connects the $x_1$-axis with an open traversing in $V_{\mathbf
j}^{(N)}$, and an open path $\gamma_2$ in $\cup_{t \in \xi
}B_t^{(N)\dagger }$ that connects the $x_2$-axis with an open
traversing in $V_{\mathbf j}^{(N)\dagger }$. Depending on whether
$\epsilon ({\mathbf j})=0 $ or $1$, we put
$$
y({\mathbf j})=\left\{
   \begin{array}{ll}
   \bigl( 3^{N+2}-\ell_h({\mathbf j}^*), 3^{N+2}-\ell_h({\mathbf j}^*)\bigr) ,
    & \hbox{ if } \epsilon ({\mathbf j})=0,\\
   \bigl( \ell_v({\mathbf j}^*), \ell_v({\mathbf j}^*) \bigr) ,
    & \hbox{ if } \epsilon ({\mathbf j})=1,
  \end{array}
  \right.
$$
and $\displaystyle Q({\mathbf j})= ( -3^{m+1}, 3^{m+1})^2
+y({\mathbf j})  $. Then by independence, the probability that the
open cluster in ${\mathcal B}$ which contains $\gamma_1$ and the
open cluster in ${\mathcal B}$ which contains $\gamma_2$ are
connected by an open path in $Q({\mathbf j})\setminus {\mathcal B}$
is not less
 than
$p^{c(m)}$, where $c(m)$ is a constant depending only on $m$. To be
more precise we can take $c(m) =8\cdot 3^{m+1}$. Since $\{
Q({\mathbf j}) ;\ |{\mathbf j}|=N-m \} $ are disjoint, the
probability that such an open connection exists for some ${\mathbf
j} \in C_{\mathbf 0}(Z)$ such that $|{\mathbf j}|=N-m$, is not less
than
$$
1- ( 1- p^{c(m)})^k.
$$
We take $k$ so large that $( 1- p^{c(m)} )^k < \varepsilon$. Then,
we have
\begin{eqnarray*}
&&P_p\left[
  \begin{array}{l}
   \hbox{ $x_1$-axis and $x_2$-axis are connected by an }\\
   \hbox{ open path in } \bigl( G_N(3,2) + (0,2\cdot 3^N) \bigr) \cup
        \bigl( G_N(2,3)+(2\cdot 3^N,0)\bigr) \\
  \end{array}
  \right]  \\
 &&\geq (1 - 2q_\varepsilon)(1-\varepsilon ).
\end{eqnarray*}
By the FKG inequality we have finally
\begin{eqnarray*}
 P_p\biggl[ C_{N, b\ell }\biggr] &\geq & P_p\bigl[ A_N(2,2) \bigr]^2
   \biggl[ 1 - \sqrt{1-P_p\bigl[ A_N(3,2)\bigr] }\biggr] \\
   && \times \biggl[ 1- \sqrt{1- P_p\bigl[ A_N(6,2)\bigr] }\biggr]
    (1-2q_\varepsilon )(1-\varepsilon )
\end{eqnarray*}
for sufficiently large $N$. Since $\varepsilon $ is arbitrarily
small, this completes the proof.

\section{Uniqueness of the critical probability}

In this section we prove the second equality in the statement (1)
and the statement
 (2) of Theorem\ref{theorem2}.
Since we have proven the equality $p_c(G^T)=p_c(S^T)$, we write
simply $p_c$ for $p_c(G^T)=p_c(S^T)$. We first claim that
\begin{equation}\label{eq,24}
\theta ({p_c})= 0.
\end{equation}
Assume that $\theta (p_c)>0$. Then by Lemma \ref{lemma2.2}, (2) and
Lemma \ref{lemma2.1.5}, (1), we have
$$
\lim_{n\rightarrow \infty }P_{p_c}[ A_n(3,1)]=1.
$$
Thus, for a $0<\theta <1$ and large $n\geq 1 $, we have
$$
P_{p_c}[A_n(3,1)]\geq 1 - 5^{-2}\theta .
$$
Then taking $\theta'\in (\theta ,1)$ and sufficiently small
$\varepsilon >0$, we have
$$
P_{p_c-\varepsilon }[A_n(3,1)] \geq 1-5^{-2}\theta',
$$
which, by the scaling argument, implies that
$$
\lim_{n\rightarrow \infty }P_{p_c-\varepsilon }[ A_n(3,1)]=1.
$$
By Lemma \ref{lemma2.4}, this means that $p_c-\varepsilon \geq p_c$,
a contradiction.

Combining (\ref{eq,24}) with Lemma \ref{lemma2.3}, we have
\begin{equation}\label{eq,25}
\liminf_{n\rightarrow \infty }P_{p_c}\left[ A_n(2,2)\right] < 1.
\end{equation}
This, together with Lemma \ref{lemma2.1} and the argument in the proof
of (2) of Lemma \ref{lemma2.1.5}, implies
that
\begin{equation}\label{eq,26}
\limsup_{n\rightarrow \infty}P_{p_c}\left[ A_n^*(18,1)\right]
  \geq \limsup_{n\rightarrow \infty }P_{p_c}[ A_n^*(6,2)] >0.
\end{equation}
\begin{lemma}\label{lemma4.1}
If $p<p_c$, then we have
\begin{eqnarray}
\liminf_{n\rightarrow \infty }P_p\left[ A_n(2,2)\right]&=&0, \label{eq,27}\\
\liminf_{n\rightarrow \infty }P_p\left[ A_n(1,1)\right]&=&0,
\label{eq,28}
\end{eqnarray}
\end{lemma}
From this lemma it is easy to obtain the final equality $p_c=
1-p_c(S^{T*})$. For this, it is sufficient to see the inequality
$p_c \leq 1- p_c(S^{T*}) $, since by Theorem \ref{theorem1}, we have
$1-p_c(S^{T*})\leq p_c$. But if $p<p_c$, by (\ref{eq,27}) and by the
scaling argument we have $P_p( A_n^*(3,1))$ converges to $1$
exponentially fast, which implies that $ 1-p\geq p_c(S^{T*})$.

The proof of Lemma \ref{lemma4.1} is essentially the same as
Kesten's original argument. Here, we sketch the proof of
(\ref{eq,27}). Let $\delta >0$ be a positive number such that
$$
\delta < \limsup_{n\rightarrow \infty }P_{p_c}\left[
A_n^*(18,1)\right]
$$
Then we can find a subsequence $\{ n_k\} $ such that
\begin{equation}\label{eq,29}
P_{p_c}\left[ A_{n_k}^*(18,1)\right] > \delta .
\end{equation}
By Russo's formula we have
\begin{equation}\label{eq,30}
\frac{d}{dp}P_p\left[ A_n(2,2)\right] = E_p\left[
N_{A_n(2,2)}\right]
\end{equation}
where $N_{A_n(2,2)}$ denotes the number of pivotal edges for
$A_n(2,2)$. Let us recall that an edge $e$ is pivotal for an event
$A$ in a configuration $\omega $ if and only if either of the
followings holds;
\begin{enumerate}
\item $\omega \in A$ and $\omega^e \not\in A$, or
\item $\omega \not\in A$ and $\omega^e \in A$,
\end{enumerate}
where
$$
\omega^e (f)=\left\{
  \begin{array}{ll}
    \omega (f) & \hbox{ if } f\not= e,\\
    1-\omega (e)& \hbox{ if } f=e.
  \end{array}
 \right.
$$
By (\ref{eq,30}), we have
$$
\frac{d}{dp}P_p\left[ A_n(2,2) \right] \geq
    E_p\left[ N_{A_n(2,2)} \
      \big\vert A_n(2,2)\right] P_p\left[ A_n(2,2)\right] .$$
Integrating this from $p$ to $p_c$, we obtain
\begin{equation}\label{eq,31}
P_{p_c}\left[ A_n(2,2)\right] \geq
   P_p\left[ A_n(2,2) \right]\exp\left\{
     \int_p^{p_c}E_q\bigl[ N_{A_n(2,2)} \ \big\vert \
              A_n(2,2) \bigr] dq \right\} .
\end{equation}
This is valid for all $n\geq 1$. We will show that for every $q\in
(p,p_c)$,
$$
E_q\left[ N_{A_{n_k +2}(2,2)} \big\vert A_{n_k+2}(2,2)\right] \geq
k\delta^5,
$$
where $\{ n_k\} $ satisfies the inequality (\ref{eq,29}). Clearly
this together with (\ref{eq,31}) proves ({\ref{eq,27}). We divide
$A_{n_k+2}$ into the sets that specify the lowest open left-right
crossing $r$ of $G_{n_k+2}(2,2)$. For a path $r$ in $G_{n_k+2}(2,2)$
connecting left side of  $G_{n_k+2}(2,2)$ with its right side, let
$$
E(r)=\left\{ \hbox{ $r$ is the lowest open left-right crossing in }
G_{n_k+2}(2,2)
    \right\} ,
$$
and
$$
E^*(r)=\left\{
\begin{array}{l}
\hbox{ there exists a closed dual path in }G_{n_k}^*(1,18) \hbox{ connecting }\\
\hbox{ the top side of } G_{n_k+2}^*(2,2) \hbox{ with a dual edge
which
 crosses } r
\end{array}
    \right\} .
$$
Then, by (\ref{eq,29}) for $0<q<p_c$ we have
$$
P_q\left[ E^*(r) \big\vert E(r) \right] > \delta ,
$$
by the FKG inequality since $E^*(r)$ is a decreasing event.
For $\omega \in E^*(r)$, let $\psi $ denote the left-most closed
dual path
 connecting the top side of $G_{n_k}^*(1,18)$ with a dual edge crossing
$r$. Further, let $e_\psi $ denote the edge in $r$ whose dual edge
$e_\psi^*$ is connected to $\psi $. Apparently $e_\psi $ is then a
pivotal edge for $A_{n_k+2}(2,2)$. Let $r_+$ denote the part of $r $
to the right of $e_\psi $, i.e. $r_+$ connects $e_\psi $ with the
right side of $G_{n_k+2}(2,2)$.

For $j=1, \ldots , k-1$, let $G_{n_j}(e_\psi )$ denote the subgraph
of $G_{n_k+2}(2,2)$ such that it is a shift of $G_{n_j+2}(1,1)$ and
it contains $e_\psi $. We write
$$ G_{n_j}(e_\psi )= G_{n_j+2}+ x_j(e_\psi )$$
so that $x_j(e_\psi )$ is the lower left corner point of
$G_{n_j}(e_\psi )$. Consider an annulus
$$ H_{n_j}= [-3^{n_j+1},3^{n_j+2}+3^{n_j+1}]^2 \setminus
           (-2\cdot 3^{n_j},3^{n_j+2}+2\cdot 3^{n_j})^2, $$
and let
$$
H_{n_j}(e_\psi )= H_{n_j}+x_j(e_\psi ).
$$
Note that $\{ H_{n_j}(e_\psi )\}_{1\leq j\leq k-1}$ are disjoint.
%
Let
$$
F(r,\psi ) = \left\{
  \begin{array}{l}
  \hbox{ $\psi $ is the left-most closed dual path connecting the top side }\\
  \hbox{ of $G_{n_k}^*(1,18)$ with a dual edge crossing $r$}
  \end{array}
     \right\}
$$
and for given $r$ and $\psi $, let
$$
C_j(r,\psi )=\left\{
  \begin{array}{l}
     \hbox{ there is a dual closed path in $H_{n_j}(e_\psi )$, located entirely }\\
    \hbox{ above $r$ and to the right of $\psi $, connecting $\psi $ with a }\\
    \hbox{ dual edge which crosses $r$}
  \end{array}
    \right\} .
$$
Then for $0<q<p_c$, by the FKG inequality we have
$$
P_q \left[ C_j(r,\psi ) \big\vert F(r, \psi )\right] > \delta^4
$$
Note that for $\omega \in F(r,\psi )\cap C_j(r,\psi )$, there exists
a pivotal edge for $A_{n_k+2}(2,2)$ in $H_{n_j}(e_\psi )$. Therefore
we have
\begin{eqnarray*}
\lefteqn{E_q\left[ N_{A_{n_k+2}(2,2)} \bigg\vert F(r,\psi ) \right] }\\
 &\geq & 1 + \sum_{j=1}^{k-1}E_q\left[ 1_{C_j(r,\psi )} \big\vert F(r,\psi ) \right]
 \\
 &\geq & 1 + (k-1)\delta^4,
\end{eqnarray*}
 since $e_\psi $ is pivotal on $F(r,\psi )$.
Therefore, we have
$$
E_q( N_{A_{n_k+2}(2,2)}\mid \hbox{ $r$ is the lowest left-right open
crossing }) \geq k\delta^5.
$$

\end{document}